\documentclass[12pt,leqno,draft]{article}

\newtheorem{theorem}{Theorem}
\newtheorem{lemma}[theorem]{Lemma}
\newtheorem{proposition}[theorem]{Proposition}
\newtheorem{definition}[theorem]{Definition}
\newtheorem{corollary}[theorem]{Corollary}

\newcommand{\begintheorem}{\addtocounter{equation}{1}\begin{theorem}}
\newcommand{\beginlemma}{\addtocounter{equation}{1}\begin{lemma}}
\newcommand{\beginproposition}{\addtocounter{equation}{1}\begin{proposition}}
\newcommand{\begindefinition}{\addtocounter{equation}{1}\begin{definition}}
\newcommand{\begincorollary}{\addtocounter{equation}{1}\begin{corollary}}

\begin{document}

\title{Notes on normed algebras, 3}

\author{Stephen William Semmes	\\
	Rice University		\\
	Houston, Texas}

\date{}

\maketitle

	Let $A$ be a countably-infinite commutative semigroup
with an identity element $0$.  Thus $A$ is equipped with a binary
operation $+$ which is commutative and associative, and
$x + 0 = x$ for all $x \in A$.

	If $f_1$, $f_2$ are two complex-valued functions on $A$, then
we can try to define the convolution of $f_1$, $f_2$ by
\begin{equation}
	(f_1 * f_2) (z) = \sum_{x + y = z} f_1(x) \, f_2(y),
\end{equation}
where more precisely the sum is taken over all $x, y \in A$ such that
$x + y = z$.  This sum makes sense if at least one of $f_1$, $f_2$ has
finite support, which is to say that it is equal to $0$ for all but at
most finitely many $x \in A$.  In particular, if we restrict our
attention to functions with finite support, then we get a nice
commutative algebra.

	For each $a \in A$ define $\delta_a$ to be the function on $A$
which is equal to $1$ at $a$ and to $0$ at all other points in $A$.
For any function $f$ on $A$, the convolution of $f$ with $\delta_0$ is
equal to $f$ again.  The convolution of $\delta_a$, $\delta_b$ is
equal to $\delta_{a + b}$ for all $a, b \in A$.  Thus the semigroup
$A$ is embedded into the convolution algebra of functions on $A$ with
finite support in such a way that the semigroup operation on $A$
corresponds exactly to convolutions of functions.  Functions on $A$
with finite support are the same as the finite linear combinations of
the $\delta_a$'s, so that convolution for any two functions on $A$
with finite support is determined by convolution of the $\delta_a$'s
and linearity.

	Consider the homomorphisms of this convolution algebra onto
the complex numbers.  In other words, these are the linear mappings
$\phi$ from the vector space of functions on $A$ with finite support
into the complex numbers such that $\phi(f_1 * f_2)$ is equal to
$\phi(f_1)$ times $\phi(f_2)$ for all functions $f_1$, $f_2$ on $A$
with finite support and $\phi(\delta_0) = 1$.  If $\phi$ is such a
homomorphism, then $\Phi(a) = \phi(a)$ defines a homomorphism from $A$
into the multiplicative semigroup of complex numbers such that
$\Phi(0) = 1$.  Conversely, if $\Phi$ is a homomorphism of $A$ into
the multiplicative semigroup of complex numbers, which means that
$\Phi(a + b) = \Phi(a) \, \Phi(b)$ for all $a, b \in A$, and if
$\Phi(0) = 1$, then we get a homomorphism $\phi$ from the convolution
algebra of functions on $A$ with finite support onto the complex
numbers by reversing the process.

	For instance, suppose that $A$ is the semigroup of nonnegative
integers.  Then for each complex number $z$ we get a homomorphism
$\Phi$ from $A$ into the multiplicative semigroup of complex numbers
such that $\Phi(0) = 1$, by setting $\Phi(z) = z^j$ when $j \ge 1$.
If $A$ is the group of integers under addition, then for each nonzero
complex number $z$ we get a homomorphism $\Phi$ from $A$ into the
multiplicative semigroup of complex numbers by setting $\Phi(z) = z^j$
for all integers $j$.

	Suppose again that $A$ is a countably infinite commutative
semigroup with identity element $0$.  A function $f$ on $A$ is said to
be summable if $\sum_{a \in A} |f(a)|$ is finite, and in this event we
write $\|f\|_1$ for the sum.  If $f_1$, $f_2$ are two functions on $A$
and at least one of $f_1$, $f_2$ is summable and the other is bounded,
then the convolution $f_1 * f_2$ can be defined as before, and is a
bounded functon on $A$.  If $f_1$, $f_2$ are two summable functions on
$A$, then the convolution $f_1 * f_2$ is also summable, and $\|f_1 *
f_2\|_1$ less than or equal to the product of $\|f_1\|_1$ and
$\|f_2\|_2$.  Thus the vector space of summable functions on $A$
becomes a commutative algebra with respect to convolution.

	Let $\Phi$ be a homomorphism from $A$ into the multiplicative
semigroup of complex numbers such that $\Phi(0) = 1$.  Assume also
that $\Phi$ is bounded, which implies that $|\Phi(a)| \le 1$ for all
$a \in A$.  Then we can define a linear mapping $\phi$ from the vector
space of summable functions on $A$ into the complex numbers by saying
that $\phi(f)$ is equal to $\sum_{a \in A} \Phi(a) \, f(a)$ for all
summable functions $f$ on $A$.  Thus $\phi(\delta_0) = 1$, and one can
check that $\phi(f_1 * f_2)$ is equal to the product of $\phi(f_1)$
and $\phi(f_2)$ when $f_1$, $f_2$ are summable functions on $A$.  We
also have that $|\phi(f)| \le \|f\|_1$ for all summable functions $f$
on $A$.

	For instance, if $A$ is the semigroup of nonnegative integers,
and $z$ is a complex number such that $|z| \le 1$, then we get such a
homomorphism $\Phi$ from $A$ into the multiplicative semigroup of
complex numbers by setting $\Phi(0) = 1$, $\Phi(j) = z^j$ when $j \ge
1$.  If $A$ is the commutative group of all integers, then for each
complex number $z$ with $|z| = 1$ we obtain such a homomorphism by
setting $\Phi(j) = z^j$ for all $j$.

	Let $A$ be a countably infinite commutative semigroup with
identity element $0$ again.  Assume that for each $a \in A$ there are
only finitely many pairs $a_1, a_2 \in A$ such that $a_1 + a_2 = a$.
In this case convolution can be defined as before for any two
functions $f_1$, $f_2$ on $A$, and the space of complex-valued
functions on $A$ becomes an algebra with respect to convolution.  For
instance, this condition holds when $A$ is the semigroup of
nonnegative integers.  However, we cannot convert a homomorphism
$\Phi$ from $A$ into the multiplicative semigroup of complex numbers
into a homomorphism $\phi$ from the convolution algebra into the
complex numbers except in trivial situations where $\Phi(a) = 0$ for
all but finitely many $a \in A$.

	One can consider continuous versions of these notions as well.
For instance, let $n$ be a positive integer, and let $A$ be a closed
convex cone in ${\bf R}^n$.  Thus $A$ is a closed subset of ${\bf
R}^n$ which contains $0$ and which has the property that $x + y \in A$
when $x, y \in A$ and $t \, x \in A$ when $x \in A$ and $t$ is a
nonnegative real number.  Let us also assume that $A$ is not contained
in a lower dimensional subspace of ${\bf R}^n$.  This is equivalent to
saying that the linear span of $A$ is all of ${\bf R}^n$, and to $A$
containing a nonempty open set.

	In general, for two functions $f_1$, $f_2$, the convolution
$(f_1 * f_2)(x)$ is defined by integrating $f_1(y)$ times $f_2(x - y)$
with respect to $y$.  For this one needs suitable conditions
on $f_1$, $f_2$, e.g., both are continuous and one has compact
support, or one is integrable and the other is bounded.  If both
$f_1$, $f_2$ have support in $A$, then the convolution does too.

	If $\zeta$ is an element of ${\bf C}^n$, then we get a
continuous homomorphism from ${\bf R}^n$ as a commutative group with
respect to addition into the nonzero complex numbers as a
multiplicative group by setting $\Phi(x)$ equal to $\exp (\zeta \cdot
x)$, where $\zeta \cdot x$ is defined to be $\sum_{j=1}^n \zeta_j \,
x_j$.  This leads to a linear mapping from continuous functions with
compact support on ${\bf R}^n$ into complex numbers by putting
$\phi(f)$ equal to the integral of $\Phi(x)$ times $f(x)$.  If $f_1$,
$f_2$ are two continuous functions with compact support on ${\bf
R}^n$, then the convolution $f_1 * f_2$ is also a continuous function
with compact support, and we get an algebra with respect to
convolution.  One can check that $\phi(f_1 * f_2) = \phi(f_1) \,
\phi(f_2)$ in this case, so that $\phi$ defines a homomorphism from
the convolution algebra of continuous functions on ${\bf R}^n$ into
the complex numbers.

	If $f_1$, $f_2$ are integrable on ${\bf R}^n$, then the
convolution $f_1 * f_2$ is defined as an integrable function on ${\bf
R}^n$, and we get a convolution algebra again.  To get a homomorphism
into the complex numbers we can start with a homomorphism $\Phi$ from
${\bf R}^n$ as an additive group into the multiplicative group of
nonzero complex numbers by setting $\Phi(x)$ equal to $\exp (i \, \eta
\cdot x)$, where now $\eta$ is an element of ${\bf R}^n$ so that
$\Phi$ is bounded, and in fact $|\Phi(x)| = 1$ for all $x \in {\bf
R}^n$.  If $f$ is integrable on ${\bf R}^n$, then we define
$\phi(f)$ to be the integral of $\Phi(x)$ times $f(x)$.  As before,
$\phi(f_1 * f_2)$ is equal to the product of $\phi(f_1)$ and $\phi(f_2)$
when $f_1$, $f_2$ are integrable, so that $\phi$ defines a homomorphism
from the convolution algebra into the complex numbers.

	Now let us restrict our attention to integrable functions
$f_1$, $f_2$ which are equal to $0$ on ${\bf R}^n \backslash A$.
Suppose that $\zeta = \xi + i \, \eta \in {\bf C}^n$, where $\xi, \eta
\in {\bf R}^n$.  Suppose also that $-\xi \in A^*$, which is to say
that $\xi \cdot x \ge 0$ for all $x \in A$.  Then $\Phi(x) = \exp
(\zeta \cdot x)$ defines a bounded homomorphism from $A$ as an
additive semigroup into the multiplicative group of nonzero complex
numbers.  One can put $\phi(f)$ equal to the integral of $\Phi(x)$
times $f(x)$ where $f$ is an integrable function supported on $A$ to
get a homomorphism from the convolution algebra of integrable
functions supported on $A$ into the complex numbers.


\begin{thebibliography}{3}


\bibitem {1} Y.~Katznelson, {\it An Introduction to Harmonic
Analysis}, second edition, Dover, 1976.

\bibitem {2} W.~Rudin, {\it Fourier Analysis on Groups},
Wiley Classics Library, 1990.

\bibitem {3} E.~Stein and G.~Weiss, {\it Introduction to Fourier
Analysis on Euclidean Spaces}, Princeton Mathematical Series {\bf 32},
Princeton University Press, 1971.




\end{thebibliography}
\end{document}